\documentclass[graybox]{SNmult}

\usepackage{type1cm}        
\usepackage{makeidx}         
\usepackage{graphicx}        
\usepackage{multicol}        
\usepackage[bottom]{footmisc}

\usepackage{newtxtext}       %
\usepackage[varvw]{newtxmath}       

\usepackage{bm}
\usepackage{subcaption}
\makeindex             

\begin{document}
	
\title*{A Comparison of Multirate Co-Simulation	Techniques for Field-Circuit Coupled Problems}
\author{Michael Wiesheu\orcidID{0000-0001-9161-5497} \\ Sebastian Schöps\orcidID{0000-0001-9150-0219} \\ Idoia Cortes Garcia\orcidID{0000-0001-6557-103X}}
\institute{Michael Wiesheu \at TU Darmstadt, Computational Electromagnetics, \email{michael.wiesheu@tu-darmstadt.de}
	\and Name of Second Author \at TU Darmstadt, Computational Electromagnetics \email{sebastian.schoeps@tu-darmstadt.de}
	\and Idoia Cortes Garcia \at Eindhoven University of Technology, Dept. of Mechanical Engineering, \email{i.cortes.garcia@tue.nl}}
%
%
\maketitle
\abstract*{This paper compares three different multirate splitting approaches for the application on field-circuit coupled magnetoquasistatic simulations. For these methods, again three different variants for exchanging values between the field and circuit are tested, namely voltages, currents and flux correction terms. All scenarios are applied on two different benchmark problems, i.e. a coil inductor and transformer model coupled to different circuits. The convergence behavior of different time steppers (Implicit Euler and Trapezoidal Rule) is determined for all possible settings, and guidelines for practical applications are derived.
}

\abstract{This paper compares three different multirate splitting approaches for the application on field-circuit coupled magnetoquasistatic simulations. For these methods, again three different variants for exchanging values between the field and circuit are tested, namely voltages, currents and flux correction terms. All scenarios are applied on two different benchmark problems, i.e. a coil inductor and transformer model coupled to different circuits. The convergence behavior of different time steppers (Implicit Euler and Trapezoidal Rule) is determined for all possible settings, and guidelines for practical applications are derived.
}

\newcommand{\mvpFEM}{\mathbf{a}}
\newcommand{\mvpFEMdot}{\dot{\mvpFEM}}

\newcommand{\elPot}{\varphi}
\newcommand{\elPotVec}{\bm{\elPot}}
\newcommand{\elPotVecDot}{\dot{\elPotVec}}

\newcommand{\current}{i} 
\newcommand{\currentVec}{\mathbf{\current}}

\newcommand{\voltage}{u} 
\newcommand{\voltageVec}{\mathbf{\voltage}}

\newcommand{\field}{\mathrm{f}}
\newcommand{\circuit}{\mathrm{c}}

\newcommand{\Mass}{M}
\newcommand{\MassMat}{\mathbf{\Mass}}
\newcommand{\Stiff}{K}
\newcommand{\StiffMat}{\mathbf{\Stiff}}
\newcommand{\rhs}{f}
\newcommand{\rhsVec}{\mathbf{\rhs}}
\newcommand{\Jac}{J}
\newcommand{\JacMat}{\mathbf{\Jac}}

\newcommand{\X}{\mathbf{X}}
\newcommand{\IncVec}{\mathbf{A}} 
\newcommand\m{\mathrm{m}}

\newcommand{\Am}{\IncVec_\m}

\newcommand{\UI}{$u$--$i$ }
\newcommand{\UPhi}{$u$--$\Phi_\mathrm{eddy}$ }
\newcommand{\IPhi}{$i$--$\Phi_\mathrm{eddy}$ }

\newcommand{\mw}[1]{\textcolor{blue}{#1}}

\newcommand{\LW}{0.5}
\newcommand{\HP}{4cm}

\newcommand{\logLogSlopeTriangle}[5]
{
	\pgfplotsextra
	{
		\pgfkeysgetvalue{/pgfplots/xmin}{\xmin}
		\pgfkeysgetvalue{/pgfplots/xmax}{\xmax}
		\pgfkeysgetvalue{/pgfplots/ymin}{\ymin}
		\pgfkeysgetvalue{/pgfplots/ymax}{\ymax}
		\pgfmathsetmacro{\xArel}{#1}
		\pgfmathsetmacro{\yArel}{#3}
		\pgfmathsetmacro{\xBrel}{#1-#2}
		\pgfmathsetmacro{\yBrel}{\yArel}
		\pgfmathsetmacro{\xCrel}{\xArel}
		\pgfmathsetmacro{\lnxB}{\xmin*(1-(#1-#2))+\xmax*(#1-#2)} 
		\pgfmathsetmacro{\lnxA}{\xmin*(1-#1)+\xmax*#1} 
		\pgfmathsetmacro{\lnyA}{\ymin*(1-#3)+\ymax*#3} 
		\pgfmathsetmacro{\lnyC}{\lnyA+#4*(\lnxA-\lnxB)}
		\pgfmathsetmacro{\yCrel}{\lnyC-\ymin)/(\ymax-\ymin)} 
		\coordinate (A) at (rel axis cs:\xArel,\yArel);
		\coordinate (B) at (rel axis cs:\xBrel,\yBrel);
		\coordinate (C) at (rel axis cs:\xCrel,\yCrel);
		\draw[#5]   (A)-- node[pos=0.5,anchor=north] {\footnotesize{1}}
		(B)-- 
		(C)-- node[pos=0.5,anchor=west] {\footnotesize{#4}}
		cycle;
	}
}

\section{Introduction}

Applying co-simulation techniques to low-frequency field-circuit coupled problems is a common way to speed up multirate simulations. Instead of a large monolithic system of equations, the individual subproblems, i.e., field and circuit, are solved independently. In the computational electromagnetics community this idea is known since the 90s, see e.g. \cite{Bedrosian_1993aa}. 
It enables the use of different tools, time discretizations, and solvers for field and circuit. The field equations, usually discretized with Finite Elements (FE), are expensive to solve. At the same time, the field is commonly assumed to be less dynamic compared to the computationally cheap but fast switched circuits, which are described in SPICE-like simulators by Modified Nodal Analysis (MNA). Decoupling through co-simulation promises a significant speed-up by treating the distinct problems separately.

The eddy-current problem with stranded conductors is denoted by
\begin{equation}
	\mathbf{F}_\field(\mvpFEM, \currentVec, \voltageVec) = 
	\begin{cases}
		\begin{aligned}
			\MassMat_\sigma \mvpFEMdot + \StiffMat_\nu \mvpFEM & = \X_\m \currentVec \\
			\X^\top_\m \mvpFEMdot &= \voltageVec\, ,
		\end{aligned}
	\end{cases}
\end{equation}
where $\MassMat_\sigma$ denotes the (singular) conductivity matrix, $\StiffMat_\nu$ the reluctivity matrix, and $\mvpFEM$ the unknowns for the magnetic vector potential \cite{Salon_1995aa}. Here, the system is a differential algebraic equation (DAE). Furthermore,  $\X_\m$ represents the coupling matrix that distributes the lumped currents $\currentVec$ in the field domain, and connects field changes with the voltage drop $\voltageVec$~\cite{Schops_2013aa}. The time integration of the second equation determines the fluxes. The DC resistance of the stranded conductor is considered in the circuit equations
\begin{equation}
	\mathbf{F}_\circuit(\elPotVec, \currentVec, \voltageVec) = 
	\begin{cases}
		\begin{aligned}
			\MassMat_\circuit \elPotVecDot + \StiffMat_\circuit \elPotVec + \Am\currentVec & = \mathbf{f}_\circuit(t) \\
			\Am^\top \elPotVec &= \voltageVec\, ,
		\end{aligned}
	\end{cases}
\end{equation}
where $\MassMat_\circuit$ denotes the (singular) mass matrix, $\StiffMat_\circuit$ the circuit stiffness matrix, $\mathbf{f}_\circuit$ circuit source terms, and $\elPotVec$ the circuit unknowns containing the electric node potentials and possible unknowns for currents across inductors and voltage sources~\cite{Gunther_2005aa}. The coupling to the field is performed through the incidence matrix $\Am$, which connects the currents to the corresponding nodes and extracts the respective voltages from the nodal potentials.

\begin{figure}[t]
	\centering
	\begin{subfigure}[c]{0.49\textwidth}
		\centering
		\includegraphics[width=\linewidth]{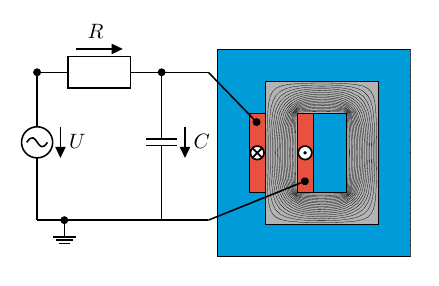}\\[-1em]
		\caption{Coupled coil inductor.}
		\label{fig:wiesheu:problems:problem1}
	\end{subfigure}%
	\hfill
	\begin{subfigure}[c]{0.50\textwidth}
		\centering
		\includegraphics[width=\linewidth]{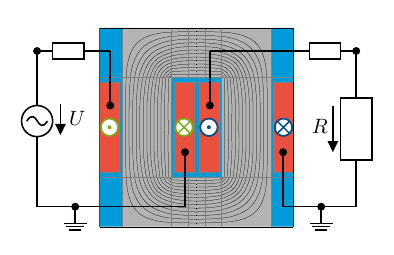}\\[-1em]
		\caption{Mutually coupled transformer.}
		\label{fig:wiesheu:problems:problem2}
	\end{subfigure}
	\caption{Field-circuit coupled benchmark problems.}
	\label{fig:wiesheu:problems}
\end{figure}

To test different co-simulation approaches, the two benchmark problems shown in Fig.~\ref{fig:wiesheu:problems} are investigated. The first problem, given in Fig.~\ref{fig:wiesheu:problems:problem1}, is a 2D inductor with a nonlinear massive iron core, connected to an $R$--$C$ circuit through only differential components. The second -- more complex -- example, given in Fig.~\ref{fig:wiesheu:problems:problem2}, is a mutually coupled 2D FE transformer with a nonlinear massive iron core connected to a circuit with no additional dynamics. As time steppers, the Implicit Euler (IE) and Trapezoidal Rule (TR) are compared since they represent widely used methods in circuit simulation. The goal is to derive rules from numerical experiments that predict which method works reliably under which circumstances. 
In the most general case, Waveform Relaxation \cite{White_1985aa} can be used to iterate until convergence is reached. This may not be ideal in terms of computational efficiency, as problems are solved repeatably, thereby reducing its potential advantage w.r.t. a monolithic approach. Furthermore, software frameworks need to allow restarts with reinitializations which is not the case for many SPICE-like solvers. We therefore concentrate on approaches that do not iterate.

The investigated methods for information exchange are explained in Sect.~\ref{sec:wiesheu:VariableExchange}.
Sect.~\ref{sec:wiesheu:CircuitField} deals with the co-simulation order `circuit-field', i.e. the circuit is solved first and its solution is passed to the field. The opposite order is explored in Sect.~\ref{sec:wiesheu:FieldCircuit}. These methods can be interpreted as multirate Lie-Trotter splittings \cite{Trotter_1959aa} or multirate co-simulation without iteration \cite{Bartel_2004ab}.
Sect.~\ref{sec:wiesheu:Strang} addresses a Strang-like splitting \cite{Strang_1968aa}. Finally, a summary and outlook is given in Sect.~\ref{sec:wiesheu:Summary}.

\section{Exchange of Variables}
\begin{figure}[t]
	\centering
	\begin{subfigure}[t]{0.35\textwidth}
		\centering
        \includegraphics[]{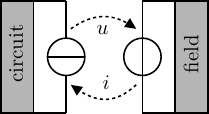}
		\caption{Voltage-current coupling.}
		\label{fig:wiesheu:cosimulation:V-I}
	\end{subfigure}%
	\hfill
	\begin{subfigure}[t]{0.32\textwidth}
		\centering
        \includegraphics[]{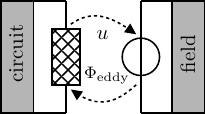}
		\caption{Voltage-flux coupling.}
		\label{fig:wiesheu:cosimulation:V-Phi}
	\end{subfigure}%
	\hfill
	\begin{subfigure}[t]{0.32\textwidth}
		\centering
        \includegraphics[]{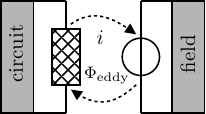}
		\caption{Current-flux coupling.}
		\label{fig:wiesheu:cosimulation:I-Phi}
	\end{subfigure}%
	\caption{Investigated coupling schemes for co-simulation.}
	\label{fig:wiesheu:cosimulation}
\end{figure}

When separating field and circuit equations, the order in which the systems are solved, and variables to be exchanged must be chosen. The investigated variants are shown in Fig.~\ref{fig:wiesheu:cosimulation}. 
\label{sec:wiesheu:VariableExchange}%
Classically, voltages and currents are exchanged, see Fig.~\ref{fig:wiesheu:cosimulation:V-I}. The field system is excited by voltages to avoid DAE problems \cite{Schops_2010ab}. This leads to
\begin{align}
	\mathbf{F}_\field(\mvpFEM, \currentVec_\field, \voltageVec_\field) &=  \mathbf{0},  &
	\voltageVec_\field &= \voltageVec_\circuit, \\
	\mathbf{F}_\circuit(\elPotVec, \currentVec_\circuit, \voltageVec_\circuit) &=  \mathbf{0},  &
	\currentVec_\circuit &= \currentVec_\field, 
\end{align}
where $\mathbf{F}_\field$ and $\mathbf{F}_\circuit$ have been defined  before. 
This approach is easy to implement, but does not exploit physical properties of the system. 
Alternatively, we prescribe a (fixed) inductance $\mathbf{L}$ to approximate the field's behavior, and provide a correction $\Phi_\mathrm{eddy}$, see Fig.~\ref{fig:wiesheu:cosimulation:V-Phi}, viz.
\begin{align}
	\Phi &= \mathbf{L}\currentVec - \Phi_\mathrm{eddy}, & 
	\Phi_\mathrm{eddy} &= \mathbf{L}\currentVec - \X^\top \mvpFEM.
\end{align}
Errors due to nonlinearities and eddy currents are compensated by updating $\Phi_\mathrm{eddy}$ accordingly. The field is still excited through a voltage and  $\Phi_\mathrm{eddy}$ is passed to the circuit. 
With a slight abuse of notation the scheme is given by 
\begin{align}
	\mathbf{F}_\field(\mvpFEM, \Phi_{\mathrm{eddy},\field}, \voltageVec_\field) &=  \mathbf{0},  &
	\voltageVec_\field &= \voltageVec_\circuit, \\
	\mathbf{F}_\circuit(\elPotVec, \Phi_{\mathrm{eddy},\circuit}, \voltageVec_\circuit) &=  \mathbf{0},  &
	\Phi_{\mathrm{eddy},\circuit} &= \Phi_{\mathrm{eddy},\field}.
\end{align}
The inductance matrix can be computed, e.g. from a linear simulation.

As third option, we use the same approach with the exception of exciting the field with currents $\currentVec$ instead of voltages, see Fig.~\ref{fig:wiesheu:cosimulation:I-Phi}. The motivation is to have a less dynamic input quantity for the field system, e.g. when comparing with pulse-width-modulated (PWM) voltages. This results in the coupling approach
\begin{align}
	\mathbf{F}_\field(\mvpFEM, \Phi_{\mathrm{eddy},\field}, \currentVec_\field) &=  \mathbf{0},  &
	\currentVec_\field &= \currentVec_\circuit, \\
	\mathbf{F}_\circuit(\elPotVec, \Phi_{\mathrm{eddy},\circuit}, \currentVec_\circuit) &=  \mathbf{0},  &
	\Phi_{\mathrm{eddy},\circuit} &= \Phi_{\mathrm{eddy},\field}.
\end{align}
Since we use fluxes as field outputs, DAE-related problems are assumed to be mitigated. However, a mathematical analysis is open. 

\begin{figure}[t!]
	\centering
	\begin{subfigure}[t]{0.48\textwidth}
		\centering
        \includegraphics[]{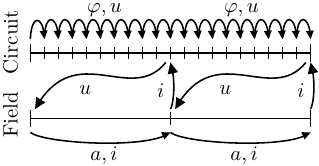}
		\caption{Circuit-Field with \UI exchange.}
		\label{fig:wiesheu:co-simulation:FieldCircuit}
	\end{subfigure}%
	\hfill
	\begin{subfigure}[t]{0.48\textwidth}
		\centering
        \includegraphics[]{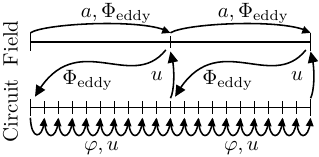}
		\caption{Field-Circuit with \UPhi exchange.}
		\label{fig:wiesheu:co-simulation:CircuiField}
	\end{subfigure}%
	\caption{Circuit-Field vs. Field-Circuit coupling with different exchange variables.}
	\label{fig:wiesheu:co-simulation:exchange}
\end{figure}

\section{Circuit-Field}
\label{sec:wiesheu:CircuitField}
The Circuit-Field variant follows the `fastest first' paradigm. This is demonstrated in Fig.~\ref{fig:wiesheu:co-simulation:FieldCircuit} for the exchange of voltages and currents. The circuit is assumed to exhibit a more dynamical behavior, and is stepped by $n$ intermediate micro-steps. Here, $n=10$ is chosen for all simulations. Unknown values are extrapolated constantly. Then the field is solved in one macro-step with linear interpolation. This drastically reduces the computational effort, since the field is more computationally demanding.

\begin{figure}[t!]
    \includegraphics[]{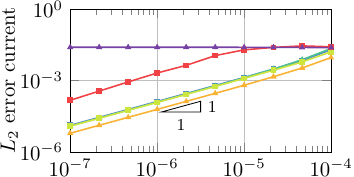}\hfill
    \includegraphics[]{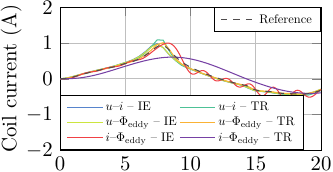}\\
    \includegraphics[]{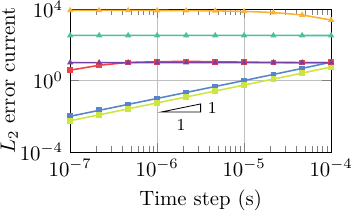}\hfill
    \includegraphics[]{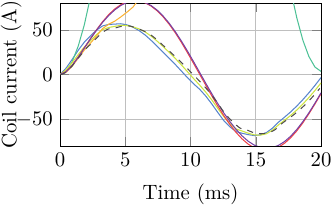}%
	\caption{Circuit-Field co-simulation convergence and coil currents calculated from the field. Top: inductor. Bottom: transformer.}
	\label{fig:wiesheu:circuitField:results}
\end{figure}

The benchmark simulations from Fig.~\ref{fig:wiesheu:problems} are excited with a sinusoidal voltage of $100$\,V at $50$\,Hz. 
The first period is simulated and compared with a very accurately resolved monolithic reference solution. The $L_2$ error of $\currentVec_\m$ in time is chosen as reference. Convergence results for the Circuit-Field simulation are found in Fig.~\ref{fig:wiesheu:circuitField:results}.

For the inductor example, the coil saturates after about 5\,ms, which causes a sudden rise in the current. The \UI exchange works for both IE and TR. However, the second order for TR is reduced to first, because the splitting is only first order accurate \cite{Hairer_2006aa}. 
The \UPhi exchange has the same convergence rates as the \UI variant. However in absolute terms, \UPhi has a lower error, presumptively because physical model information is included in its information exchange. 
The \IPhi exchange converges for IE, but very small time steps are necessary to reach an acceptable error. The correction of the current seems to be delayed, such that it oscillates around the true value. The TR case does not converge.

For the transformer example, the only methods that converge are \UI and the \UPhi variant with IE. Also the \IPhi exchange works eventually with IE, but again only for extremely small time steps. TR fails in all cases, because the circuit is purely algebraic. We conjecture that errors due to inconsistent initial conditions are accumulated. 

\section{Field-Circuit}
\label{sec:wiesheu:FieldCircuit}
The Field-Circuit variant follows the `slowest first' paradigm. The simulation starts with a macro step, where the field is computed. After that, the circuit is simulated with $n$ micro-steps, with the interpolated input provided by the field. As before we use linear interpolation and constant extrapolation. Two exemplary time windows for the \UPhi exchange are shown in  Fig.~\ref{fig:wiesheu:co-simulation:CircuiField}.

\begin{figure}[t!]
    \includegraphics[]{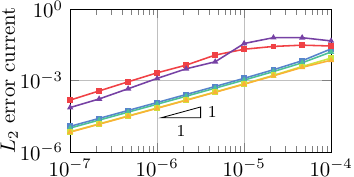}\hfill
    \includegraphics[]{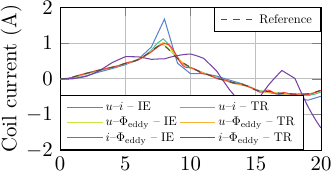}\\
    \includegraphics[]{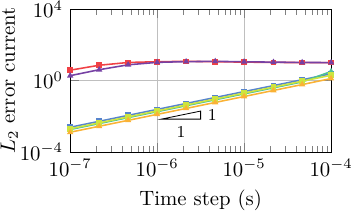}\hfill
    \includegraphics[]{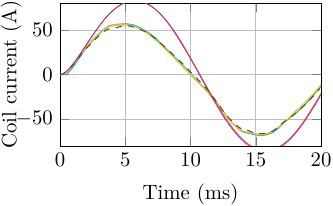}
	\caption{Field-Circuit co-simulation convergence and coil currents calculated from the circuit. Top: inductor. Bottom: transformer.}
	\label{fig:wiesheu:fieldCircuit:results}
\end{figure}

The convergence results for the two benchmark problems are shown in Fig.~\ref{fig:wiesheu:fieldCircuit:results}. For the inductor problem, \UI converges for both IE and TR, again with order one. 
%
The \UPhi exchange converges for IE as before. 
The \IPhi exchange suffers from the same problems as above, i.e. a very slow convergence is observed. However, also TR converges.

For the transformer problem, the convergence is comparable to Circuit-Field with one major difference: The \UI and \UPhi exchanges with TR converges with first order. This indicates that the initial conditions for the circuit are consistent if the field is stepped first, and the circuit may interpolate the field quantities. 

\section{Strang-Splitting}
\label{sec:wiesheu:Strang}
\begin{figure}[b!]
	\centering
    \includegraphics{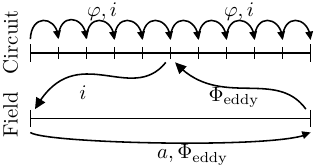}
	\caption{Strang coupling with \IPhi exchange.}
	\label{fig:wiesheu:co-simulation:Strang}
\end{figure}

Strang splitting alternates between a half step for the first subproblem, a full step for the second one, and another half step for the first subproblem \cite{Strang_1968aa}. This yields a splitting order of two for ordinary differential equations. Our systems, however, are DAEs. In the context of field-circuit coupling and multirate problems, a Strang-like splitting could be implemented as shown in Fig.~\ref{fig:wiesheu:co-simulation:Strang} for one time window with the \IPhi exchange. 
The circuit is started with $n/2$ micro-steps. The end values are passed to the field, which performs one macro-step. Finally, the time window is completed by another circuit simulation with $n/2$ micro-steps with input values from the field. Here, we do not interpolate, but use the last value to stay close to the Strang setting.

\begin{figure}[t!]
    \includegraphics[]{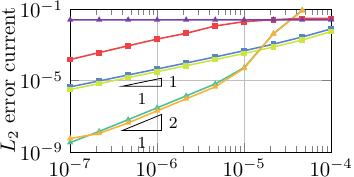}\hfill
    \includegraphics[]{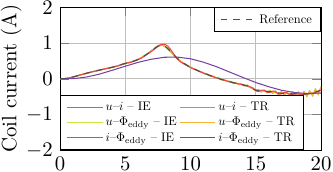}\\
    \includegraphics[]{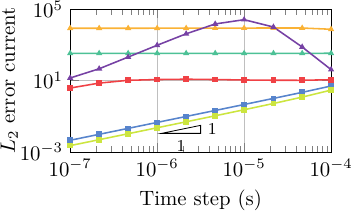}\hfill
    \includegraphics[]{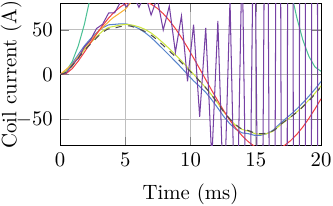}
	\caption{Strang co-simulation convergence and coil currents calculated from the field. Top: inductor. Bottom: transformer.}
	\label{fig:wiesheu:strang:results}
\end{figure}

Applying this method on the benchmark problems yields the convergence results presented in Fig.~\ref{fig:wiesheu:strang:results}. For the inductor problem, the \UI exchange has optimal convergence for both IE and TR, i.e. order 1 and 2, respectively. Even though the problem is not an ODE, the differential components connected to the exchange variables, i.e., the capacitor in parallel, seem to conserve the second order convergence property of the Strang-splitting. However, oscillations occur for larger time steps, where IE outperforms the TR. 
The same conclusions are reached for the \UPhi exchange, which displays the same convergence behavior, but tends to perform better in absolute terms.
The \IPhi exchange again only converges using IE with small time steps and not at all for TR.

For the transformer problem, the \UI exchange converges for IE. Applying TR does not converge and again inconsistent initial values are assumed to be the reason. The same conclusion is reached for the \UPhi exchange. Again, the \IPhi exchange does not converge with TR, and small time steps are necessary for IE to converge.

\section{Summary and Outlook}
\label{sec:wiesheu:Summary}

Three different multirate co-simulation methods were applied to two field-circuit coupled problems. Eddy currents have been included in the field formulation, however, magnetostatic field models behave analogously.

It has been shown that the \UI exchange, where voltages are fed to the magnetic field and currents are given to the electric circuit, is the most robust method. IE works reliably in all scenarios, TR converges in some cases with (its optimal) second order. The \UPhi exchange, where voltages are fed to the field, and a flux correction term is passed back to the circuit alongside a fixed inductance, also converges in all tested cases with IE, and in certain scenarios with TR. Compared to \UI, the \UPhi variant has lower absolute errors, which can be explained by the additional information encoded through the inductance. A third exchange method has been proposed, namely the \IPhi variant, where currents are passed to the field, and a flux correction term is passed to the circuit. The motivation for this method is to have a continuous field input, since voltages might be discontinuous, e.g. in the case of PWM signals. While the \IPhi exchange converges using IE, very small time steps are necessary to reach satisfactory errors. Convergence with TR could not always be reached for this case. In general, applying TR is more challenging because of oscillations which we associate with inconsistent initial conditions.

Future work will focus on deepening the mathematical understanding for the observed convergence behavior, and testing further methods that allow for the multirate co-simulation of PWM signals.

	\begin{acknowledgement}
		The work  is supported by the joint DFG/FWF Collaborative Research Centre CREATOR (CRC -- TRR361/F90) at TU Darmstadt, TU Graz and JKU Linz as well as the Graduate School CE within the Centre for Computational Engineering at TU Darmstadt.
	\end{acknowledgement}
	\ethics{Competing Interests}{
		The authors have no conflicts of interest to declare that are relevant to the content of this chapter.}

\end{document}